\documentclass[a4paper,11pt]{article}
\usepackage[para]{bigfoot}
\usepackage[utf8]{inputenc}
\usepackage{eurosym,multicol,times}
\usepackage{color}
\usepackage{url}
\usepackage{amsmath,amssymb}
\usepackage{hyperref}
\usepackage{listings}
\usepackage{graphicx}

\newcommand{\Sectiontitle}[1]{\bigskip\bigskip\vspace*{-6mm}%
\noindent\begin{minipage}[t]{\columnwidth}\Large\bfseries\hrule\vspace{1.5mm}%
\begin{center}{#1}\end{center}\vspace{-2.8mm}\hrule\end{minipage}\medskip}

\newcommand{\Subsectiontitle}[1]{\ \\ \noindent{\bfseries #1}\\ \noindent}

\DeclareNewFootnote[para]{default}

\begin{document}
\setcounter{footnote}{0}
\setcounter{figure}{0}

  \begin{center}
    \begin{minipage}[b]{115mm}
      {\Large\bfseries FindStat --\\[5px] \hspace*{10px} the combinatorial statistics database}  \\[1em]
      {\bfseries 
        Chris Berg (Google)\\
        Viviane Pons (Universit\"at Wien)\\
        Travis Scrimshaw (UC Davis)\\
        Jessica Striker (North Dakota State University)\\
        Christian Stump (Freie Universit\"at Berlin)      
      }  \\[1em]
      \href{mailto:info@findstat.org}{\texttt{info@findstat.org}}
    \end{minipage}
    \begin{minipage}[b]{57mm}
    \includegraphics[height=45mm]{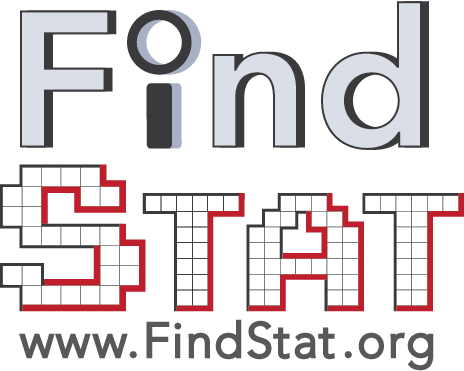}
    \end{minipage}
  \end{center}

\vspace{-10pt}

\begin{multicols}{2}
\noindent

\thispagestyle{empty}
\pagestyle{empty}

The FindStat project provides an online platform for mathematicians, particularly for combinatorialists, to gather information about combinatorial statistics and their relations.  
As of January 2014, the FindStat database contains 173 statistics on 17 combinatorial collections.
The project was initiated by Chris Berg and Christian Stump in 2011 at the Laboratoire de combinatoire et d'informatique mathématique, Université du Québec à Montréal, Canada.
In 2013, Viviane Pons, Travis Scrimshaw, and Jessica Striker joined the project.

\Sectiontitle{An example}

Combinatorial statistics arise naturally all over mathematics. Before we give the definition of a combinatorial statistic, let us first describe an example. Suppose you were studying a collection of polynomials $\{f_\pi\}$ indexed by permutations, and the degree of the first few polynomials were as follows:
\begin{gather*}
  deg(f_{12}) = 0 \quad
  deg(f_{21}) = 1 \\[5pt]
  deg(f_{123}) = 0 \quad 
  deg(f_{132}) = 2 \quad 
  deg(f_{213}) = 1 \\[0pt]
  deg(f_{231}) = 2 \quad 
  deg(f_{312}) = 1 \quad 
  deg(f_{321}) = 3 \\[5pt]
  deg(f_{1234}) = 0 \quad 
  deg(f_{1243}) = 3 \quad
  deg(f_{1324}) = 2 \\[0pt]
\end{gather*}
\vspace{-25pt}

You would like to know what underlying rule determines the degree of the polynomial. The savvy combinatorialist will tell you that the degree of $f_\pi$ appears to be the \emph{major index} of $\pi$,
$$maj(\pi) = \sum_{1 \leq i < n \atop \pi_i > \pi_{i+1}} i.$$

What's that? You say you aren't a savvy combinatorialist? Never heard of a Mahonian statistic and are now riddled with the shame and dismay of your peers? Relax! This is what FindStat was built for! Just point your browser to
\begin{center}
  \small \url{www.FindStat.org/StatisticFinder}
\end{center}
and simply input the degree of each of your polynomials. Wait a few seconds, and FindStat will tell you that you're looking at the major index of a permutation.

\smallskip

But what if you know that your polynomials are indexed by permutations, but you cannot figure out which polynomial belongs to which permutation of a given length?
Even in this case, FindStat will be able to tell that you are looking for a statistic that has the same distribution (or generating function) as the major index.

\Sectiontitle{Online Databases}

\Subsectiontitle{The OEIS} We know what you're thinking! ``Wait a second, isn't this exactly what the OEIS does?''
Indeed, the success of the \emph{On-Line Encyclopedia of Integer Sequences} (OEIS)~\cite{OEIS} was a primary motivation in the creation of FindStat.
The OEIS database contains more than $200,000$ integer sequences
 and has had an immense  impact on the work of mathematicians in all fields of mathematics.
There were two main reasons for starting the FindStat project, the first one being more of an annoyance, the second being more profound:

\begin{itemize}
  \item {\bf Combinatorial collections do not necessarily admit a canonical total ordering} which would be needed to index a given combinatorial statistic in the OEIS. In the example above, we study permutations. While there are many ways to index permutations (e.g. lexicographically), there is no default. Even more annoying, some combinatorial objects, like graphs, have much less obvious total orderings.
  
  \item {\bf Combinatorial collections usually come with much more structure} than only a linear ordering. In particular, there are many natural relationships between various combinatorial collections. A statistic on one combinatorial collection can often be transformed to a meaningful statistic on an entirely different collection of objects. In fact, those following along online may have  noticed that FindStat didn't exclusively return the major index of a permutation as a solution to the above problem.
As of January 1, 2014, searching the FindStat database for this data yields about $22,000$ database searches after following more than $1,200$ combinatorial maps. The database results provide the major index\footnote{\url{www.FindStat.org/St000004}}, as expected, as well as connections to the number of inversions\footnote{\url{www.FindStat.org/St000018}} and multiple other statistics in the database.
\end{itemize}

\Subsectiontitle{Fingerprint databases for theorems}
As discussed by S.C.~Billey and B.E.~Tenner in~\cite{BT2013}, FindStat is an example of a
\begin{center}
  \emph{fingerprint database for theorems.}
\end{center}
Roughly speaking, this means that FindStat is actually providing a canonical way of storing and recognizing theorems; when I search for one statistic and find another, FindStat is implying that these statistics are in fact equal or related through well-described connections. Moreover, similar to a single sequence from the OEIS, a statistic page often contains information regarding the contexts in which the statistic arises, and references containing further information.

\Sectiontitle{Project outline}

A \emph{combinatorial collection} $\mathcal{S}$ is considered to be a set with interesting combinatorial properties, a \emph{combinatorial map} is a map between two combinatorial collections that has an explicit combinatorial description, and a \emph{combinatorial statistic} on $\mathcal{S}$ is a map $st$ that associates an integer to each element in $\mathcal{S}$, $st:\mathcal{S}\longrightarrow\mathbb{Z}$.

\Subsectiontitle{The main aims} The two main aims of the FindStat project are to

\medskip

\noindent 1. provide an online platform to gather information about combinatorial collections and statistics and their relations. This includes
\begin{itemize}
  \item adding new combinatorial statistics to the FindStat database, and
  \item filling the corresponding wiki with information about combinatorial collections and statistics.
\end{itemize}

\noindent 2. provide a web interface to
\vspace*{-5pt}
\begin{itemize}
  \item test if your data is a known combinatorial statistics in the database,
  \item test if your data can be obtained from known combinatorial statistics in the database by applying combinatorial maps, or
  \item test if your data has the same distribution as a known combinatorial statistic in the database by applying combinatorial maps.
\end{itemize}
\vspace*{-5pt}
\Subsectiontitle{Generating functions} The first question many people asked about FindStat was whether it also supports generating functions -- and now FindStat does!
We have already seen in the introductory example that you can search for exact values attached to the given permutations. 
The idea is that even if you have only the \emph{distribution} of the degrees of your polynomials among all permutations of a given size, 
you can still search for statistics with the same generating function
$$F_n(q) = \sum_{\pi \in \operatorname{Perm}(n)} q^{deg(f_\pi)}.$$
FindStat not only supports generating functions, but also has a very flexible user interface which allows you to search
for statistics matching the distribution of your input data according to any partition of your data.

\Subsectiontitle{Contributing to the database} It might have happened in the introductory example that the degrees of the polynomials were not found in the database.
You should then think of adding your data through the web interface
\begin{center}
  \small\url{www.FindStat.org/NewStatistic}
\end{center}
so that the next person searching for this statistic will be aware of the context in which you found it.

We tried to make it as simple as possible to contribute. All you need to provide is your data, the context in which you encountered the statistic, and, if possible, references.
The most important part is to provide the correct data; we recommend using your favorite computer algebra system (or the Sage cell\footnote{\url{https://github.com/sagemath/sagecell}} provided by the web interface) to produce the data for you.
If you do not have an algorithm computing the data for you, you can alternatively provide a few values by hand. Maybe the next person to find your statistic in the database will be able to provide further values.

\Subsectiontitle{The FindStat wiki} In addition to the pure database of combinatorial statistics, the FindStat project serves as a platform to provide information about combinatorial statistics and maps. For example, information about permutation statistics and combinatorial maps for permutations can be found at \small\url{www.FindStat.org/Permutations} \normalsize and its subsequent pages.
Anyone can contribute to the wiki to improve the content and the presentation of the information.


\Subsectiontitle{License and technical information} Contributions to the FindStat project, i.e., to the database and to the wiki, are licensed under a \emph{Creative Commons License}\footnote{\href{http://creativecommons.org/licenses/by-sa/3.0/deed.en_US}{Creative Commons Attribution-ShareAlike 3.0 Unported License}.}.
The main technologies used in this project are a \emph{MoinMoin wiki engine} for serving the website~\cite{MOIN}, and a server running \emph{Sage} in which the combinatorial maps are implemented and the actual computations are performed~\cite{Sage}.


\thispagestyle{empty}

\end{multicols}
\end{document}